\documentclass[12pt]{article}

\usepackage{amsfonts,amssymb,amsmath,amsthm,eucal}

\newcommand{\la}{\lambda}
\newcommand{\Ln}{\Lambda_n}
\newcommand{\de}{\delta}
\newcommand{\Dc}{\mathcal{D}_n}

\newcommand{\lang}{\left\langle}
\newcommand{\rang}{\right\rangle}
\newcommand{\lb}{\left[}
\newcommand{\rb}{\right]}
\newcommand{\wh}{\widehat}

\newcommand{\lt}{\vartriangleleft}

\newcommand{\gt}{\vartriangleright}

\newcommand{\mm}[2]{\mod\{#1\}_{#2}}

\newtheorem{Proposition}{Proposition}
\newtheorem{Theorem}{Theorem}
\newtheorem{Corollary}{Corollary}

\begin{document}

\title{A remark on Fourier pairing and binomial formula 
for Macdonald polynomials} 
\author{Andrei Okounkov\thanks{
 Department of Mathematics, University of California at
Berkeley, Evans Hall \#3840, 
Berkeley, CA 94720-3840. E-mail: okounkov@math.berkeley.edu}
}
\date{} 
\maketitle

\begin{abstract}
  We give a concise direct proof of the orthogonality of 
interpolation Macdonald polynomials 
with respect to the Fourier pairing and  briefly
discuss some immediate applications of this orthogonality, such
as the symmetry of the Fourier pairing and the binomial
formula. 
\end{abstract}

\section{Introduction}

Fourier pairing, introduced by Cherednik \cite{C1}, is a
fundamental notion in the theory of Macdonald polynomials. 
In its simplest instance, it pairs the algebra $\Ln$ of
symmetric polynomials in $n$ variables with the algebra
$\Dc$ of Macdonald commuting difference operators  acting
on $\Ln$ \cite{M1}. By definition, 
\begin{equation}
\lang D, f\rang = \lb D \cdot f \rb (\wh 0) \,, \quad D\in\Dc,f\in\Ln\,,
\label{<>0}
\end{equation}
where $\wh 0$ is a certain distinguished point. There is 
a natural isomorphism $\Dc\cong\Ln$, which makes \eqref{<>0} a
quadratic form on $\Ln$. The most important and useful property
of this form is its symmetry, see \cite{C1,M1,M2}. 

The main observation of this note is that there is a very 
natural orthogonal basis for the form \eqref{<>0}.
Namely, this is the basis $\{I_\mu\}$ of  
the interpolation Macdonald polynomials, which have
been intensively
studied by Knop, Olshanski, Sahi, the author, and others,
see for example \cite{K1,O1,S1} and also \cite{O4}, 
and references therein. The
polynomials $I_\mu$ are defined by some very simple multivariate
Newton-type interpolation conditions and have found many
remarkable applications. 

The orthogonality of $I_\mu$ with respect to \eqref{<>0},
stated as Theorem \ref{T1} below, 
follows rather easily from the definitions, without using
any nontrivial properties of the polynomials $I_\mu$. 
It implies at
once that \eqref{<>0} is symmetric. 

Also, the orthogonality of $I_\mu$ gives
immediately the expansion of the simultaneous eigenfunctions 
of $\Dc$, known as the symmetric Macdonald polynomials $P_\la$,
in the basis $\{I_\mu\}$. This
expansion, reproduced in Theorem \ref{T3} below, 
is the binomial formula for $P_\la$, see \cite{O2}. In fact, the 
orthogonality of $I_\mu$ is essentially equivalent to 
the binomial theorem, but it certainly appears to be
a much more basic, natural, and appealing property.

The binomial theorem of \cite{O2} has been extended to 
a more general setting, including other classical 
root systems and the nonsymmetric Macdonald polynomials,
see \cite{K2,O3,S2}. We are not trying to pursue the greatest
possible generality in this note. Our intention, rather, 
is to show how the basic idea works in the simplest
nontrivial example of the usual symmetric Macdonald
polynomials. We even consider the (almost) trivial
one-dimensional case first to give a completely
elementary illustration of what is going on.

It should be pointed out that there is another 
source of orthogonality relations for the 
polynomials $I_\mu$. Namely, the polynomials
$I_\mu$ can be obtained from the symmetric
Macdonald polynomials of type $BC_n$ (this 
can be seen very explicitly by degenerating
the binomial formula of \cite{O3} to the 
binomial formula for the polynomials $I_\mu$
\cite{O2}). 

This note is based my unpublished paper written in 
Fall of 1997. Later, it was intended to be a part
of a survey article on interpolation Macdonald
polynomials on which we were working together 
with G.~Olshanski. I would like to thank NSF
(grant DMS-0096246), Sloan foundation, and Packard
foundation for financial support.

\section{Simplest example}

\subsection{}

As a warm-up, let us consider the one-dimensional case first.
Let the operator $T$ act on polynomials in $x$ by
$$
\lb T f \rb(x) = f(qx) \,.
$$ 
Obviously, the monomials $x^n$, $n=0,1,2,\dots$, are the
eigenfunctions of this operator with eigenvalues
$q^n$\,. Consider the following bilinear form
\begin{equation}
\lang g, f\rang = \lb g(T)\cdot f \rb (1)  \,.
\label{<>1}
\end{equation}
In the basis $\{x^n\}$, this form has the matrix
$$
\big[ q^{nm} \big]_{n,m=0,1,\dots} = 
\begin{bmatrix}
  1 & 1 & 1 & 1 & \hdots \\
 1 & q & q^2 & q^3 & \hdots \\
 1 & q^2 & q^4 & q^6 & \hdots \\
1 & q^3  & q^6 & q^9 & \hdots \\
\vdots & \vdots &\vdots & \vdots & \ddots 
\end{bmatrix}
$$
and is clearly symmetric. 

It is also obvious that 
\begin{equation}
x^* = T \,, \quad T^* = x \,,
\label{F1}
\end{equation}
where $x$ denotes the operator of multiplication by the 
independent variable and
star denotes the adjoint operator with respect to \eqref{<>1}.
Clearly, \eqref{F1} is an anti-automorphism of 
the $q$-Heisenberg algebra generated by $T$ and $x$ 
(subject to the relation $Tx=qxT$) and deserves to be called the
\emph{Fourier transform}.

\subsection{} 

Now consider the polynomial 
\begin{equation}
I_n = (x-1)(x-q) \cdots (x-q^{n-1})\,, \quad n=0,1,\dots \,, \label{In}
\end{equation}
satisfying the following Newton interpolation conditions
\begin{align}
 &I_n  \equiv  x^n \mod \{x^m\}_{m<n} \,, \label{I1} \\
 &I_n(q^m) = 0\,,  \quad 0\le m < n \,. \label{I2}
\end{align}
We have the following

\begin{Proposition} The polynomials $I_n$ are orthogonal 
with respect to the form \eqref{<>1}, namely
\begin{equation}
  \label{o1}
 \lang I_n, I_m \rang = \de_{n,m} \, I_n(q^n) \,.
\end{equation}
\end{Proposition}

\begin{proof}
We will intentionally avoid using the symmetry of \eqref{<>1}
in our argument because, in the general case, we want to 
obtain the analogous symmetry as a corollary. 

It is clear from definition \eqref{<>1} that 
\begin{equation}
\lang x^n, f\rang = f(q^n)  \,,
\label{xf}
\end{equation}
and since $T \cdot x^n = q^n \, x^n$ we also have
$$
\lang g, x^n\rang = g(q^n)  \,. 
$$
{}From \eqref{I2} it now follows that 
$$
\lang x^m, I_n \rang = \lang I_n, x^m \rang  = 0\,,  \quad m < n \,,
$$
and it is also clear that 
$$
\lang x^n, I_n \rang = \lang I_n, x^n \rang  = I_n(q^n)\,.
$$
Now the property \eqref{I1} concludes the proof.
\end{proof}

The following expansion is immediate from \eqref{o1} and \eqref{xf}
$$
x^n = \sum_m \frac{\lang x^n, I_m\rang} 
{\lang I_m, I_m\rang} \, I_m(x) =
\sum_m \frac{I_m(q^n) \, I_m(x)}{I_m(q^m)} \,.
$$
This is the Newton interpolation of $x^n$ with nodes
$1,q,q^2,\dots$ and also an instance of the $q$-binomial
theorem.

\section{Symmetric Macdonald polynomials} 

\subsection{} 

We now turn to polynomials in $n$ variables $x_1,\dots,x_n$. 
We denote 
$$
\lb T_i f \rb(x_1,\dots,x_n) = f(x_1,\dots,qx_i,\dots,x_n) \,.
$$
Let $t$ be an additional parameter and introduce, following
Macdonald \cite{M1}, the following operators
$$
D_k = t^{k(k-1)/2} 
\sum_{|S|=k} \, d_S(x) \, \prod_{i\in S} T_i \,,
$$
where the summation is over subsets $S\subset\{1,\dots,n\}$
of cardinality $k$ and 
$$
d_S(x) = \prod_{i \in S, \, j\notin S} 
\frac{t x_i - x_j}{x_i - x_j}  \,. 
$$

\subsection{} 

The operators $D_k$ commute and take symmetric polynomial to symmetric
polynomials. They act triangularly in the basis of 
monomial symmetric functions, namely, 
$$
D_k \cdot m_\la  \equiv  e_k(\wh\la) \, m_\la \mm{m_\mu}{\mu<\la} \,,
$$
where $\la=(\la_1\ge \la_2 \ge \dots \ge \la_n)$ is
a partition,
$$
m_\la = x_1^{\la_1} \cdots x_n^{\la_n} \quad  + 
\quad\textup{permutations} 
$$
is the corresponding monomial symmetric function, 
$e_k=m_{(1^k)}$ is the
$k$th elementary symmetric function, $\wh\la$ denotes
the following point
$$
\wh\la = \left(q^{\la_1} t^{n-1}, q^{\la_2} t^{n-2}, \dots,
q^{\la_{n-1}} t, q^{\la_n}\right) \,,
$$
and  $\mu < \la$ stands for the dominance order on partitions
$$
\mu \le  \la \Leftrightarrow 
\left(
  \begin{matrix}
    \mu_1 \le \la_1\,, \\
\mu_1 + \mu_2 \le \la_1 + \la_2\,,\\ \dots \\
\mu_1 +\dots +\mu_n = \la_1 +\dots +\la_n 
  \end{matrix}
\right) \,. 
$$
The simultaneous 
eigenfunctions of $D_k$ 
\begin{equation}
D_k \cdot P_\la = e_k(\wh\la) \, P_\la \,,
\label{DP}
\end{equation}
normalized by
\begin{equation}
P_\la \equiv  m_\la  \mm{m_\mu}{\mu<\la} 
\label{Pm}
\end{equation}
are known as the Macdonald symmetric polynomials.

\subsection{}

Let $\Ln$ denote the algebra of symmetric polynomials in 
$n$ variables. It is clear from \eqref{DP} that the 
map 
$$
D: \Ln \owns e_k \mapsto D_k 
$$
extends to an algebra homomorphism such that
\begin{equation}
D(g)\cdot P_\la = g(\wh\la) \, P_\la \,, \quad g\in \Ln\,. 
\label{DgP}
\end{equation}
We now define, following Cherednik \cite{C1}, the following 
\emph{Fourier pairing}
\begin{equation}
  \label{<>}
  \lang g, f\rang = \lb D(g)\cdot f \rb (\wh0)\,,
\quad f,g\in\Ln  \,.
\end{equation}
This is an analog of \eqref{<>1}. It is clear that
$$
\lang h g, f\rang = \lang g, D(h)f\rang \,.
$$
In other words, $D(h)=h^*$, where $h$ is considered as a multiplication
operator and star denotes its \emph{Fourier transform}, that is,
the  adjoint operator with respect to \eqref{<>}. It is also 
clear from 
\eqref{DP}  that the pairing \eqref{<>}  takes
the normalized eigenfunction 
\begin{equation}
N_\la = \frac{P_\la}{P_\la(\wh 0)}
\label{N}
\end{equation}
to the $\delta$-function at $\wh\la$, namely 
\begin{equation}
  \label{gP}
  \lang g, N_\la \rang = g(\wh\la) \,.
\end{equation}

\subsection{} 

Our goal is now to produce an explicit orthogonal basis 
for the quadratic form \eqref{<>}. As in \eqref{In}, this basis
will consist of certain Newton interpolation polynomials.

Let $\lt$ be any total ordering of the set of partitions $\la$
compatible with both the ordering of partitions by their size
$|\la|$ and, for partitions of the same number, the dominance
ordering. Define the interpolation Macdonald polynomials 
$I_\mu$ by the following generalization  of \eqref{I1} and 
\eqref{I2}
\begin{align}
 &I_\mu  \equiv m_\mu \mm{m_\la}{\la\lt\mu} \,, \label{Im1} \\
 &I_\mu(\wh\la) = 0\,,  \quad \la \lt \mu \,. \label{Im2}
\end{align}
The existence and uniqueness of such polynomials for generic $q$ and $t$ is 
clear from their existence and uniqueness for $t=1$,
which is elementary.

\subsection{} 

It can be shown, see for example \cite{K1,O1,S1} and also 
\cite{O4}, that the polynomials
$I_\mu$ do not depend on the choice of the ordering $\lt$ 
and satisfy the much stronger \emph{extra vanishing} property
\begin{equation}
  \label{exv}
  I_\mu(\wh\la) = 0\,,  \quad \mu\not\subset\la  \,.
\end{equation}
By the binomial formula \eqref{bin}, this gives the 
following strengthening of \eqref{Im1} 
\begin{equation}
  \label{Ism1}
  I_\mu  \equiv P_\mu \mm{P_\la}{\la\subset \mu} \,. 
\end{equation}
The extra vanishing \eqref{exv} will not, however, be needed for 
what follows making our argument applicable in the 
situations where the analog of \eqref{exv} is not available. 

\subsection{} Our main result is the following 

\begin{Theorem}\label{T1} The polynomials $I_\mu$ are orthogonal with
respect to the Fourier pairing \eqref{<>}. 
\end{Theorem}

An immediate corollary of this theorem is the following
central result of the theory of Macdonald polynomials

\begin{Corollary}[Koornwinder,\cite{M1}]\label{C1}
The Fourier pairing \eqref{<>} is symmetric. 
\end{Corollary}

Koornwinder actually proved an equivalent symmetry, namely the
following label-argument symmetry for the normalized 
polynomials \eqref{N}
$$
N_\la (\wh \mu) = N_\mu (\wh \la) \,.
$$
Numerous application of this symmetry, such as, for example, 
Pieri-type formulas for Macdonald polynomials, can be found in 
\cite{C1,M1,M2}. 

\subsection{}\label{pf1}
The proof of Theorem \ref{T1} goes in two steps. First, we 
claim that
$$
\lang I_\mu , I_\la \rang = 0 \,, \quad \mu \gt \la \,.
$$
Indeed, by \eqref{Im1}, \eqref{Pm}, and \eqref{gP}, this is equivalent to 
$$
\lang I_\mu , N_\la \rang = I_\mu(\wh\la)=0  \,, \quad \mu \gt \la \,,
$$
which is indeed true by \eqref{Im2}. 

\subsection{}\label{pf2}
Now we prove that 
$$
\lang I_\mu , I_\la \rang = 0 \,, \quad \mu \lt \la \,.
$$
By \eqref{Im1} this is equivalent to proving that
$\lang m_\mu , I_\la \rang = 0$ if $\mu \lt \la$. 
Since
$$
e_\mu \overset{\textup{\it def}}= e_{\mu_1} \cdots e_{\mu_n}
\equiv m_\mu \mm{m_\nu}{\nu<\mu} \,,
$$
it suffices to prove that
$$
\lang e_\mu , I_\la \rang = 0 \,, \quad \mu \lt \la \,.
$$
By definition \eqref{<>}, this is equivalent to 
\begin{equation}
\lb D_\mu \cdot I_\la \rb (\wh 0) = 0 \,, 
\quad D_\mu= D_{\mu_1} \cdots D_{\mu_n} \,,\label{DI}
\end{equation}
which will now be established. 

\subsection{} 

It is a crucial property of the operators $D_k$ that
$$
\left(
  \begin{array}{l}
   \la_i=\la_{i+1}\,, \\
   i\notin S,\, i+1\in S
  \end{array}
\right) 
\Rightarrow d_S(\wh \la) = 0 \,.
$$
It follows that 
$$
\lb D_k \cdot f \rb (\wh \la)  = 
\sum_{\nu/\la = \textup{vertical $k$-strip}}
d_{S(\nu,\la)}(\wh \la) \, f(\wh \nu)  \,, 
$$
where $S(\nu,\la) = \{ i, \nu_i > \la_i\}$\,. 
It follows that 
\begin{equation}
\lb D_\mu \cdot f \rb (\wh 0) = 
\sum_{\nu \le \mu} c_{\mu,\nu} \, f(\wh \nu)\,,
\label{Dmf}
\end{equation}
for some coefficients $c_{\mu,\nu}$. A similar property can be
established in more general context, such as e.g.\ for nonsymmetric
Macdonald polynomials \cite{C2}. 

It is clear that \eqref{Dmf} together with \eqref{Im2} imply
\eqref{DI} and this concludes the proof of Theorem \ref{T1}. 

\subsection{}
Theorem \ref{T1} can be sharpened as follows

\begin{Theorem}\label{T2} We have
  \begin{equation}
    \label{eT2}
    \lang I_\mu, I_\nu \rang =
\de_{\mu,\nu} \, I_\mu(\wh \mu) \, P_\mu(\wh 0)  = \de_{\mu,\nu} \,
c_{\mu,\mu} \, I_\mu(\wh \mu) \,. 
  \end{equation}
\end{Theorem}

In particular, this shows that $P_\mu(\wh 0)=c_{\mu,\mu}$ which,
after making the number $c_{\mu,\mu}$ explicit, can
be seen to be equivalent to the known formula for $P_\mu(\wh 0)$, see
\cite{M1}. 

\begin{proof}
  Arguing as in Section \ref{pf1}, we see that
$$
\lang I_\mu, I_\mu \rang = \lang I_\mu, P_\mu \rang = 
I_\mu(\wh \mu) \, P_\mu(\wh 0) \,.
$$
On the other hand, arguing as in Section \ref{pf2}, we get
$$
\lang I_\mu, I_\mu \rang = 
\lb D_\mu \cdot I_\mu \rb (\wh 0) = c_{\mu,\mu} \, I_\mu(\wh \mu) \,.
$$
\end{proof}

\subsection{} 
Theorem \ref{T2} implies the following Newton interpolation 
formula
$$
f = \sum_\mu \frac{\lang I_\mu, f\rang}{\lang I_\mu, I_\mu\rang} \, I_\mu \,,
\quad f\in\Ln \,. 
$$
In particular, applying this to $N_\la$ and using \eqref{gP}
we obtain the following expansion (in which we explicitly 
kept the variables $x$ in order to stress the label-argument
symmetry). 

\begin{Theorem}[Binomial theorem, \cite{O2}]\label{T3} We have
\begin{equation}
N_\la(x)  = \sum_{\mu} \frac{I_\mu(\wh \la) \, I_\mu(x)}
{\lang I_\mu, I_\mu\rang} \,. 
\label{bin}
\end{equation}
\end{Theorem}

It follows from \eqref{exv} that only those $\mu$ such that
$\mu\subset\la$ actually appear in this expansion.

\end{document}